\documentclass[12pt]{amsart}
\usepackage{amsmath,amsthm,amsfonts,amssymb,eucal}

\renewcommand {\a}{ \alpha }
\renewcommand{\b}{\beta}

\newcommand{\g}{\gamma}

\renewcommand{\d}{\delta}
\newcommand{\s}{\sigma}

\renewcommand{\L}{\Lambda}
\newcommand{\z}{\zeta}

\newcommand{\p}{\partial}
\newcommand{\om}{\omega}
\newcommand{\Om}{\Omega}

\newcommand{\oq}{\ {\raise 7pt\hbox{${\scriptstyle\circ}$}}
	\kern -7pt{
		\hbox{$Q$}}}

\newcommand{\R}{ \mathbb R}

\newcommand{\Rd}{ \mathbb R^d}

\newcommand {\GA}{\mathfrak A}
\newcommand {\GD}{\mathfrak D}

\newcommand {\GG}{\mathfrak G}
\newcommand {\GH}{\mathfrak H}

\newcommand {\GS}{\mathfrak S}

\newcommand {\GW}{\mathfrak W}

\newcommand {\bx}{\mathbf x}

\newcommand {\bk}{\mathbf k}

\newcommand {\by}{\mathbf y}

\newcommand {\bn}{\mathbf n}

\newcommand {\bxi}{\boldsymbol\xi}

\newcommand{\CD}{\mathcal D}

\newcommand{\plainW}[2]{\textup{{\textsf{W}}}^{#1, #2}}
\newcommand{\plainC}[1]{\textup{{\textsf{C}}}^{#1}}

\newcommand{\plainS}{\textup{{\textsf{S}}}}

\newcommand{\plainL}[1]{\textup{{\textsf{L}}}^{#1}}

\DeclareMathOperator{\tr}{{tr}}

\newcommand{\1}
{{\,\vrule depth3pt height9pt}{\vrule depth3pt height9pt}
	{\vrule depth3pt height9pt}{\vrule depth3pt height9pt}\,}

\DeclareMathOperator {\im }{{Im}}
\DeclareMathOperator {\re} {{Re}}

\DeclareMathOperator{\op}{{Op}}

\newtheorem{thm}{Theorem}[section]
\newtheorem{cor}[thm]{Corollary}
\newtheorem{lem}[thm]{Lemma}
\newtheorem{prop}[thm]{Proposition}
\newtheorem{cond}[thm]{Condition}

\theoremstyle{definition}
\newtheorem{defn}[thm]{Definition}

\newtheorem{rem}[thm]{Remark}

\numberwithin{equation}{section}

\newcommand{\bee}{\begin{equation}}
	\newcommand{\ene}{\end{equation}}
\newcommand{\bees}{\begin{equation*}}
	\newcommand{\enes}{\end{equation*}}
\newcommand{\bes}{\begin{split}}
	\newcommand{\ens}{\end{split}}

\newcommand{\bet}{\begin{thm}}
	\newcommand{\ent}{\end{thm}}
\newcommand{\bel}{\begin{lem}}
	\newcommand{\enl}{\end{lem}}
\newcommand{\bec}{\begin{cor}}
	\newcommand{\enc}{\end{cor}}
\newcommand{\bep}{\begin{proof}}
	\newcommand{\enp}{\end{proof}}
\newcommand{\ber}{\begin{rem}}
	\newcommand{\enr}{\end{rem}}

\begin{document}
\hoffset -4pc

\title
[Functions of self-adjoint operators]
{{Functions of self-adjoint operators in ideals of compact operators}}
 \author{Alexander V. Sobolev}
 \address{Department of Mathematics\\ University College London\\
Gower Street\\ London\\ WC1E 6BT UK}
\email{a.sobolev@ucl.ac.uk}
\keywords{Compact operators, quasi-normed ideals, functional calculus,  
	Wiener-Hopf operators with discontinuous symbols, 
	quasi-classical asymptotics}
	\subjclass[2010]{Primary  47G30,  47B15; Secondary 45M05, 35S05, 47B10, 47B35}
	
\begin{abstract}
For self-adjoint operators $A, B$, a bounded operator $J$, and 
a function $f:\R\to\mathbb C$ we obtain bounds in 
quasi-normed ideals of compact operators for the difference $f(A)J-Jf(B)$ 
in terms of the operator  $AJ-JB$. 
The focus is on 
functions $f$ that are smooth everywhere except for finitely many points. 
 A typical example is the function $f(t) = |t|^\gamma$ 
 with $\gamma \in (0, 1)$. 
 The obtained results are applied to 
 derive a two-term quasi-classical asymptotic formula for the trace 
 $\tr f(S)$ with $S$ being a Wiener-Hopf operator with a discontinuous symbol. 
\end{abstract}

\maketitle

\section{Introduction} 
In this paper we study a pair of self-adjoint operators $A$, $B$ on  Hilbert spaces $\GH$ and $\GG$ 
respectively. 
We are interested in  estimates in various quasi-normed ideals of compact operators 
for the ``quasi-commutators"  of the form $f(A) J - J f(B)$ in terms of the 
``perturbation"  $AJ-JB$, where $J: \GG\to\GH$ is a bounded operator 
and $f:\R\to \mathbb C$ is a  suitable function. 
There is a vast literature concerned with problems of this type, with 
a large number of deep results. Our intention is to improve 
some of the existing estimates 
for a very specific class of functions $f$. 
The focus will be on continuous functions $f$ that are smooth everywhere except 
possibly for finitely many points. 
One example of such function is $f(t) = |t|^\gamma$ with $\g >0$.  
In this introduction we do not provide a detailed survey of the 
known results but 
concentrate on the directly relevant ones only, 
further references can be found e.g. in \cite{AP} and \cite{PS}. 
By $\GS$ we denote a (quasi)-normed two-sided ideal of compact operators, 
and by $\GS_p, 0<p<\infty$ -- the classical Schatten-von Neumann ideals. 

In \cite{Pel1, Pel3} it was found 
that if $f$ belongs to the 
Besov class $B^1_{\infty 1}(\R)$ then the function $f$ is 
$\GS_1$-\textit{operator-Lipschitz}, i.e.
\begin{equation}\label{peltrace:eq}
\| f(A) - f(B)\|_{\GS_1}\le C\|A-B\|_{\GS_1}, C=C(f),
\end{equation}
for arbitrary self-adjoint operators $A, B$ such that $A-B\in\GS_1$. 
Conversely, as shown in \cite{Pel1}, the 
estimate \eqref{peltrace:eq} implies that 
$f\in B^1_{1 1}(\R)$ locally. 
Paper \cite{FP} identifies a meaningful class of self-adjoint operators, for which
the condition  $f\in B^1_{1 1}(\R)$ is also sufficient for 
\eqref{peltrace:eq}.

For the Schatten-von Neumann classes 
$\GS_p$, $1<p<\infty$ conditions on the function $f$ look simpler. Precisely,  
for arbitrary uniformly Lipschitz functions $f$ 
it was shown in \cite{PS} that 
\begin{equation}\label{PS:eq}
\| f(A)-f(B)\|_{\GS_p}\le c_p \|f\|_{\textup{\tiny Lip}}
\| A-B\|_{\GS_p}, \ \ 
\|f\|_{\textup{\tiny Lip}} = \sup_{x\not = y} \frac{|f(x) - f(y)|}{|x-y|}.
\end{equation}
The classes $\GS_p$ with $p\in (0, 1)$ 
were studied in \cite{Pel} for unitary operators $A$ and $B$. 
We discuss this in more detail in Remark \ref{unitary:rem}. 
     
The function $f(t) = |t|^\gamma$, $\g\in (0, 1)$ was studied 
in \cite{BKS}. 
Let $\GS$ be a  
normed ideal with the \textit{majorization property}, see \cite{GK} 
for the definition. 
This assumption is not too restrictive as any separable 
ideal (e.g. $\GS_p, 0<p<\infty$) possesses this property. 
As shown in \cite{BKS} (see also 
\cite{BS_DOI}), for $\g \in (0, 1)$, 
if $A\ge 0$ and $B\ge 0$ are such that $|A-B|^\g \in \GS$, 
then 
 \begin{equation}\label{BKS:eq}
 \|A^\g - B^\g\|_{\GS}\le \| |A-B|^\g\|_{\GS}.
 \end{equation}
Observe that the function $|t|^\g, \g\in (0, 1)$ belongs to 
the the  H\"older-Zygmund class $\L_\g(\R) 
= B^\g_{\infty, \infty}(\R)$ locally. 
Among other functional spaces, this space  
was considered in the recent article \cite{AP}.
In fact, \cite{AP} brings us closer to 
the objects studied in the current paper as it contains results on the 
quasi-commutators $f(A)J - Jf(B)$. 
Precisely, for any function $f\in \L_\g(\R), \g\in (0, 1)$ 
it was shown in \cite{AP} that  
\begin{equation}\label{AP:eq}
\| |f(A) J - J f(B)|^{\frac{1}{\g}}\|_\GS
\le C(f) \|J\|^{\frac{1-\g}{\g}}\|AJ-JB\|_\GS,
\end{equation}
under the assumption that the \textit{Boyd index} $\b(\GS)$ of the quasi-normed 
ideal $\GS$ is strictly less than $1$, see \cite{AP}, Theorem 11.5. 
The definition of the Boyd index can be found 
e.g. in \cite{AP}, Section 3. 
For the Schatten-von Neumann ideals $\GS_p, 0<p<\infty,$ the index is found by the simple formula $\b(\GS_p) = p^{-1}$. 

None of the results quoted above generalizes the others but   
some of them have non-empty intersections. 
Let us compare, for instance \eqref{BKS:eq} and \eqref{AP:eq} 
for the Schatten - von Neumann classes. Then \eqref{BKS:eq} 
gives  
\begin{equation}\label{BKS1:eq}
\|A^\g - B^\g\|_{\GS_p}\le \| A-B\|_{\GS_{p\g}}^\g,
\end{equation}
for any $p\ge 1$ and $\gamma\in (0, 1)$, and \eqref{AP:eq} gives (see \cite{AP}, Theorem 11.7) 
\begin{equation}\label{AP1:eq}
\| f(A)J - Jf(B)\|_{\GS_p}
\le C(f)^\g \|J\|^{1-\g} \|AJ-JB\|_{\GS_{p\g}}^\g,
\end{equation}
under the condition $p\g >1$. 
On the one hand \eqref{AP1:eq} is valid for the entire class $\L_\g(\R)$, and it 
allows $J\not = I$, but on the other hand, \eqref{AP1:eq} holds under the more restrictive assumption 
$p\g >1$.
  
One aim of this paper 
is to derive the following ``hybrid" of \eqref{BKS:eq} and \eqref{AP:eq}. 
For the sake of discussion 
we state the result in a somewhat simplified form, see 
Theorem \ref{quasinorm_est:thm} for the precise statement. 
Assume that $f\in\plainC\infty(\R\setminus\{z\}), z\in\R,$ is a compactly supported 
function satisfying the condition
\begin{equation}\label{fbound_pre:eq}
|f^{(k)}(t)|\le C_k |t-z|^{\g-k},\ k = 0, 1, \dots, \ t\not = 0,
\end{equation}
with some $\g >0$. Let $\GS$ be a quasi-normed ideal. Then for any 
$\s\in (0, \g), \s\le 1,$ the bound holds
\begin{equation}\label{main_pre:eq}
\||f(A)J-Jf(B)|^{\frac{1}{\s}}\|_\GS\le C(f)\|J\|^{\frac{1-\s}{\s}} \|AJ-JB\|_\GS.
\end{equation} 
Emphasize that in contrast to \eqref{AP:eq}, 
the value $\s = \g$ is not allowed. On the other hand, there 
are no restrictions on 
the ideal $\GS$.  

If $\g >1$ then in the formula \eqref{main_pre:eq} one can take $\s=1$. Thus  
for $\GS = \GS_p, 1<p<\infty$ and $J=I$ the bound \eqref{main_pre:eq} 
is in agreement with \eqref{PS:eq}. 
For $p\in (0, 1)$ and $J=I$ the bound \eqref{main_pre:eq} is in line with 
the results of \cite{Pel}, see Remark \ref{unitary:rem} for details. 

Since our choice of the function $f$ is very specific, 
the proof of \eqref{main_pre:eq} does not require sophisticated 
methods employed in \cite{AP, Pel1, Pel, Pel3, PS} where 
various general functional classes were studied. In particular, 
we do not make use of the Double Operator Integrals techniques. 
Instead we rely on the representation 
of $f(A)$ for a self-adjoint operator $A$ in terms 
of the \textit{quasi-analytic 
extension} of the function $f$, 
which has become known as 
the Helffer-Sj\"ostrand formula, see \cite{HelSjo, EBD}. 
The convenient quasi-analytic extension is constructed in Lemma 
\ref{ab:lem}. 


In Theorem \ref{Szego1:thm} we focus on the following useful special case 
of the bound \eqref{main_pre:eq}. 
Let $A$ be a self-adjoint operator and let $P$ be an orthogonal projection. 
Then, using \eqref{main_pre:eq} with $J=P$, $B=PAP$ we obtain the bound
\begin{equation}\label{LSg:eq}
\|Pf(PAP)P - Pf(A)P\|_\GS\le C(f) \||P A(I-P)|^\s\|_\GS.
\end{equation}
A bound of a similar nature was previously 
derived in \cite{LS1, LS2} for 
arbitrary $f\in \plainW{2}{\infty}_{\textup{\tiny loc}}(\R)$:
\begin{equation*} 
\bigl|\tr\bigl(Pf(PAP)P - Pf(A)P\bigr)\bigr|
\le \frac{1}{2} \|f''\|_{\plainL\infty} \| P A(I-P)\|_{\GS_2}^2.
\end{equation*} 
The above two inequalities are helpful in problems 
involving Szeg\H o-type  estimates and/or 
asymptotics, see e.g. \cite{LS1, LS2, LRS, Sob}.  

The last section of the paper, Sect. 4, illustrates the practical use 
of the bound \eqref{LSg:eq} with an example of a multi-dimensional 
Wiener-Hopf operator with a discontinuous symbol. The operator in question, denoted by 
$S_\a$, is defined by \eqref{WH:eq}, 
where $\a\ge 1$ is the ``quasi-classical parameter". 
The objective is to obtain a two-term asymptotics 
of the trace $\tr g(S_\a)$ with a non-smooth function $g$, as $\a\to\infty$. 
The function $g$ is allowed to have finitely many singularities of the type described 
by \eqref{fbound_pre:eq}. For smooth $g$ the sought two-term asymptotic formula 
was justified in \cite{Sob} and \cite{Sob2}.
The main result of Sect. 4 is contained in Theorem \ref{Szego3:thm}, and its proof 
consists in ``closing" the asymptotic formula for smooth $g$ with the help of the bound 
\eqref{LSg:eq}. 
 The generalization to non-smooth functions 
 is motivated, in part, by applications in information theory and statistical physics, 
 see e.g.   \cite{GiKl}, \cite{HLS}, \cite{LeSpSo}. Further discussion is deferred until Sect. 4. 

\textbf{Acknowledgements.} The author is grateful 
to W. Spitzer for useful remarks. 
This work was supported by EPSRC grant EP/J016829/1. 

\section{Main results}\label{est:sect}

\subsection{Quasi-normed ideals of compact operators} 
We need some information from the theory of ideals of compact 
operators. Details can be found 
in \cite{GK}, \cite{BS}, \cite{Pie}. 
Let $\GS\subset \GS_\infty$ be a two-sided ideal.  
Recall that a functional $\| \ \cdot\ \|_{\GS}$ 
defined for $T\in\GS$ is said to be a quasi-norm if 
\begin{enumerate}
\item
$\|T\|_\GS >0$ if $T\not = 0$,
\item
$\| z T\|_\GS = |z|\|T\|_{\GS}$ for any $z\in \mathbb C$,
\item
there exists a number $\varkappa\ge 1$ such that 
\begin{equation*}
\|T_1+T_2\|_\GS\le \varkappa \bigl(\|T_1\|_\GS + \|T_2\|_\GS\bigr).
\end{equation*}
If, in addition, the conditions below are satisfied
\item 
$\|X T Y\|_\GS\le \|X\| \ \|Y\|\ \|T\|_\GS $,  
for any bounded $X, Y$ and $A\in \GS$,
\item
$\|T\|_\GS = \|T\|$ for any one-dimensional operator $T$,
\end{enumerate}
then the quasi-norm $\| \ \cdot\ \|_\GS$ is said to be symmetric. 
The ideal $\GS$ is said to be a quasi-normed ideal if it is endowed with 
a (symmetric) quasi-norm, and is complete. 
We usually omit the term ``symmetric" for brevity. 
If $\varkappa=1$ then the quasi-norm becomes a norm.

Note an important property of quasi-norms. 
Below by $s_k(T), k = 1, 2, \dots,$ we denote 
singular numbers of the operator $T\in \GS_{\infty}$.   

\begin{lem} \label{quasi-norm:lem}
Let $T\in \GS$ and let $S\in \GS_\infty$ be operators such that 
$s_k(S)\le M s_k(T), k = 1, 2, \dots,$ with some constant $M>0$. Then $S\in \GS$ 
and $\|S\|_\GS\le M \|T\|_{\GS}$.
\end{lem}

For normed ideals this lemma was 
proved in \cite{GK}, and the proof for quasi-normed ideals 
is the same. It shows that the 
quasi-norm $\|T\|_\GS$ depends only on the singular numbers of the 
operator $T\in \GS$. 
This means in particular that $\|T\|_{\GS} = \|T^*\|_{\GS} = \||T|\|_{\GS}$, 
where $|T| = \sqrt{T^*T}$. 

We say that a quasi-normed ideal $\GS$ is a 
\textit{$q$-normed ideal} if there exists an equivalent 
quasi-norm $\|\ \cdot\ \|_\GS$ which satisfies the 
\textit{$q$-triangle inequality}: 
\begin{equation}\label{qtriangle:eq} 
\|T_1+T_2\|_{\GS}^q\le \|T_1\|_{\GS}^q + \|T_2\|_{\GS}^q,
\end{equation}
for any $T_1, T_2\in \GS$, see e.g. \cite{Pie}.  
In fact, any quasi-normed ideal $\GS$ is a 
$q$-normed ideal with the $q \in (0, 1]$ 
found from the equation $\varkappa = 2^{q^{-1}-1}$ 
($q=1$ refers to a normed ideal).

As an example, we can take as 
$\GS$ any Schatten-von Neumann ideal $\GS_p, \ p\in (0, \infty)$ with 
the standard (quasi)-norm 
\begin{equation*}
\|T\|_{\GS_p} = \biggl[\sum_{k=1}^\infty s_k(T)^p\biggr]^{\frac{1}{p}}.
\end{equation*}
If $p\ge 1$, then this functional is a norm, 
and if $p\in (0, 1)$ then it is a $p$-norm, 
see \cite{Rot} and also \cite{BS}.

\subsection{The estimates}
Let $A$ and $B$ be two self-adjoint operators acting on 
the Hilbert spaces $\GH$ and $\GG$ respectively, and 
let $J: \GG\to \GH$ be a bounded operator. Consider the form
\begin{equation*}
V[u, w] = (Ju, Aw) - ( JB u, w), u\in D(B), w\in D(A).
\end{equation*} 
Suppose that 
\begin{equation*}
|V[u, w]|\le C\|u\| \ \|w\|, 
\end{equation*}
i.e. this form defines an operator $V: D(B)\to\GH$ which 
extends to a bounded operator on the entire space $\GG$. This implies that 
$J$ maps $D(B)$ into $D(A)$. We use the notation $V = AJ-JB$. 
Let $R(z; A) = (A-z)^{-1}$, $\im z\not = 0$. Under the assumption that 
$V:\GG\to\GH$ is a 
bounded operator, we can write the resolvent identity
\begin{equation}\label{res:eq}
R(z; A) J - J R(z; B) = - R(z; A) V R(z; B).  
\end{equation} 
 
\begin{lem}\label{sandwich:lem}
Suppose that the operator $V=AJ-JB$ is such that 
$|V|^\s\in \GS$ with some $\s\in (0, 1]$. Then for all $y = \im z\not =0$ we have 
\begin{equation*}
\|R(z; A) V R(z; B)\|_{\GS}\le \||V|^\s\|_{\GS} \|J\|^{1-\s}
\frac{2^{1-\s}}{|y|^{1+\s}}.
\end{equation*}
\end{lem} 

\begin{proof}
Denote
\begin{equation*}
W=R(z; A) V R(z; B).
\end{equation*}
By definition of the quasi-norm,
\begin{equation}\label{normout:eq}
\|W\|_\GS\le \||W|^{1-\s}\| \ \||W|^\s\|_\GS = \|W\|^{1-\s}\  \||W|^\s\|_{\GS}
\le \frac{2^{1-\s}}{|y|^{1-\s}}\|J\|^{1-\s}\||W|^\s\|_{\GS},
\end{equation}
where we have used the trivial bound for the left-hand side of \eqref{res:eq}:
$\|W\|\le 2|y|^{-1}\|J\|$. 
In order to estimate the quasi-norm on the right-hand side of \eqref{normout:eq} 
estimate the singular values $s_k(|W|^\s)$: 
\begin{equation*}
s_k(|W|^\s) = s_k(W)^\s
\le |y|^{-2\s} s_k(V)^\s.
\end{equation*}
Therefore by Lemma \ref{quasi-norm:lem}
\begin{equation*}
\||W|^\s\|_\GS\le |y|^{-2\s}\||V|^\s\|_\GS.
\end{equation*}
Substituting this bound into \eqref{normout:eq} 
we get the required estimate. 
\end{proof}

We are interested in bounds for the difference 
\begin{equation*}
f(A)J - J f(B),
\end{equation*} 
where $f$ is a function satisfying the following condition.  
Below we denote by $\chi_R$ the characteristic function of 
the interval $(-R, R)$, $R>0$. 
 
\begin{cond}\label{f:cond}
Assume that for some integer $n \ge 1$ the function 
$f\in\plainC{n}(\R\setminus\{ x_0 \})\cap~\plainC{}(\R)$, $x_0\in\R$, satisfies the 
bound 
\begin{equation}\label{fnorm:eq}
\1 f\1_n = 
\max_{0\le k\le n}\sup_{x\not = x_0} |f^{(k)}(x)| |x-x_0|^{-\g+k}<\infty
\end{equation}
 with some $\g >  0$, 
and is supported on the interval $[x_0-R, x_0+R]$ with some $R>0$.  
\end{cond}

For a function $f$ satisfying the above condition 
the bound holds: 
\begin{equation}\label{fbound:eq}
|f^{(k)}(x)| \le \1 f\1_n |x-x_0|^{\g-k}\chi_R(x-x_0),  
k = 0, 1, \dots,  n, \ \quad  \   x\not = x_0.
\end{equation}  
One can immediately deduce from \eqref{fbound:eq} that 
\begin{equation}\label{H:eq}
	\begin{cases}
\|g\|_{\plainL\infty}\le \1 g\1_0 R^\g,\\[0.2cm] 
	\ \ |g(t_1)-g(t_2)|\le C_\g \1 g\1_1 |t_1-t_2|^\varkappa, \ 
	\varkappa = \min\{\gamma, 1\}, 
\end{cases}
\end{equation}
for any $t_1, t_2\in \R$,
so that $g\in \plainC{0, \varkappa}(\R)$. 
 Here by $\plainC{0, \varkappa}(\R^n), n\ge 1,$ 
 we denote the standard class of H\"older-continuous functions $f$ 
 with the finite norm
 \begin{equation*}
 	\sup_{\bx} |f(\bx)| 
 	+ \sup_{\bx\not = \by} \frac{|f(\bx) - f(\by)|}{|\bx-\by|^\varkappa}.
 \end{equation*} 
 The next theorem constitutes the main result of the paper.

\begin{thm}\label{quasinorm_est:thm} 
Suppose that $f$ 
satisfies Condition \ref{f:cond} with some $\g >0$,  $n\ge 2$ and $R >0$. 
Let $\GS$ be a $q$-normed ideal where $(n-\s)^{-1} < q\le 1$ with some 
number $\s \in (0, 1]$, $\s < \g$.

Let $A, B$ be two self-adjoint operators as described above 
such that $V = AJ-JB$ is a bounded operator.
 Suppose that 
$|V|^\s\in \GS$. Then 
\begin{equation}\label{fest:eq}
\|f(A) J - J f(B)\|_{\GS} \le C_n R^{\gamma-\s}\1 f\1_n\|J\|^{1-\s}\||V|^\s\|_{\GS},
\end{equation}
with a positive constant $C_n$ independent of the operators $A, B, J$, 
function $f$ and parameter $R$.  
\end{thm}
  
One should observe that the parameters $n$ and $\s$ in Theorem 
\ref{quasinorm_est:thm} are not entirely independent. Indeed, the condition 
$(n-\s)^{-1} < q\le 1$ does not allow $n = 2$ and $\s = 1$ at the same time.  
We'll need to remember this fact in the proof of Theorem \ref{Szego3:thm} later on. 

\begin{rem}\label{unitary:rem}
It is appropriate to compare Theorem \ref{quasinorm_est:thm} with the results of the paper \cite{Pel} mentioned in the introduction.  
In \cite{Pel} it was shown for a pair of unitary operators $U_1$ and $U_2$ 
that 
\begin{equation}\label{unitary:eq}
	f(U_1)-f(U_2)\in \GS_p \ \ 
	\textup{under the assumption that}\ \ U_1 - U_2\in\GS_p, p\in (0, 1),
\end{equation}
if $f\in B^{\frac{1}{p}}_{\infty p}(\mathbb T_1)$, where 
$\mathbb T_1$ is the unit circle. Conversely, 
\eqref{unitary:eq} implies that $f\in B^{\frac{1}{p}}_{p p}(\mathbb T_1)$. 
These conditions can certainly be appropriately 
rephrased for self-adjoint  operators with 
the help of the Cayley transform. 

For the sake of comparison, in Theorem \ref{quasinorm_est:thm}  
assume for simplicity that $f\in \plainC\infty(\R\setminus\{0\})$ is a function 
such that 
$f(t) = |t|^\g, \g>0,$ for all $|t|\le 1$ and $f(t) = 0$ for $|t|\ge 2$.
Then using \eqref{fest:eq} with $\GS = \GS_p$, $J = I$, $R=2$ and $\gamma>1$ we get that 
\begin{equation*}
\|f(A) - f(B)\|_{\GS_p}\le C_p \|A-B\|_{\GS_p},
\end{equation*}
for arbitrary $p\in (0, 1]$ 
with a constant $C_p$ independent of $A, B$. 
The chosen function $f$  belongs to 
 $B^{\nu}_{rq}(\R)$, $r\in (0, \infty], q\in (0, \infty)$,  
 if and only if   $\nu < \g+ r^{-1}$. 
 Thus 
 $f$ does not satisfy the sufficient condition 
 $f\in B^{\frac{1}{p}}_{\infty p}(\R)$ from \cite{Pel}, if $\g \le p^{-1}$. 
 On the other hand, the necessary condition 
 $f\in B^{\frac{1}{p}}_{p p}(\R)$ is satisfied for any $\g >0$.
 
\end{rem}

Note the following scaling property:

\begin{rem}\label{fest:rem}
Theorem \ref{quasinorm_est:thm} for arbitrary $R>0$ 
follows from Theorem \ref{quasinorm_est:thm} for $R=1$. 
Indeed, without loss of generality one may assume that $x_0 = 0$. 
Note that the function $g(t) = R^{-\g} f(Rt)$ 
satisfies \eqref{fbound:eq} with $R=1$ and that 
$\1 g\1_n = \1 f\1_n$. Now use 
bound \eqref{fest:eq} for the function $g$ and the operators 
$A' = R^{-1}A, B' = R^{-1} B$. 
\end{rem} 
 
It is also convenient to have a separately stated result for smooth functions $f$.

\begin{cor}\label{smooth:cor}
Suppose that $g\in \plainC{n}_0(-\rho, \rho)$, with some $\rho >0$ and $n\ge 2$. 
Let $\GS$ be a $q$-normed ideal where $(n-\s)^{-1} < q\le 1$ with some 
number $\s \in (0, 1]$.
Let $A, B$ be two self-adjoint operators as in Theorem 
\ref{quasinorm_est:thm}. Then 
	\begin{equation}\label{smooth:eq}
	\|g(A)J - J g(B)\|_{\GS}
	\le C \max_{0\le k\le n} \bigl(\rho^k \| g^{(k)}\|_{\plainL\infty}\bigr)   
	\rho^{-\s} \|J\|^{1-\s} \|||V|^\s\|_{\GS},
	\end{equation}
	with a constant $C$ independent of the operators $A, B, J$, function $g$ and 
parameter $R$. 
\end{cor}

\begin{proof} 
	Suppose first that $\rho=1$ and without loss of generality set 
	\begin{equation*}
	\max_{0\le k\le n}\|g^{(k)}\|_{\plainL\infty} =  1. 
	\end{equation*}
 Then the  function 
	\begin{equation*}
		f(t) = (t-2)^2 \bigl( g(t) (t-2)^{-2}\bigr)
	\end{equation*}
	clearly satisfies \eqref{fbound:eq} 
	with $\g = 2$, $x_0 = 2$, $R=3$ and 
	$\1 f\1_n\le C$. 
	Therefore by Theorem \ref{Szego1:thm},
	\begin{equation*}
		\|g(A)J - J g(B)\|_{\GS}
		\le C \|J\|^{1-\s}\||V|^\s\|_{\GS},
	\end{equation*}
	which proves \eqref{smooth:eq} for $\rho=1$. 
	
	If $\rho >0$ is arbitrary, then use the first part of the proof 
	for the function $f(t) = g(\rho t)$ and operators $A' = A \rho^{-1}$, 
	$B' = B \rho^{-1}$.  
\end{proof}

Now we use Corollary \ref{smooth:cor} 
to obtain bounds similar to \eqref{fest:eq} for functions with unbounded supports. 
We concentrate on smooth functions $g$ satisfying the bound
\begin{equation}\label{g:eq}
|g^{(k)}(x)|\le (1+|x|)^{-\b}, \b >0,\ k = 1, 2, \dots, n, x\in \R.
\end{equation}

\begin{cor}
	Suppose that $g$ satisfies \eqref{g:eq} with some $n \ge 2$ and $\b >0$. 
	Let the ideal $\GS$ and operators $A,B$ be as in Corollary \ref{smooth:cor}. 
	If $q\b>1$, then 
	\begin{equation}\label{gest:eq}
	\|g(A) J - J g(B)\|_{\GS} \le C_n  \|J\|^{1-\s}\||V|^\s\|_{\GS},
	\end{equation}
	with a positive constant $C_n$ independent of the 
	operators $A, B, J$ and function $g$. 
\end{cor}

\begin{proof} 
	Let $\Upsilon\in \plainC\infty_0(\R)$ 
	be a function such that 
	$\Upsilon(t) = 0$ for all $|t|\ge 1$. 
	We pick $\Upsilon$ in such a way that 
	$\sum_{m\in\mathbb Z} \Upsilon(x-m) = 1$, $x\in\R$. Let 
	$g_m(x) = \Upsilon(x-m) g(x)$, $m\in\mathbb Z$. Since
	$\|g_m^{(k)}\|_{\plainL\infty}\le C (1+|m|)^{-\b}$, 
	$k = 0, 1, \dots, n$, it follows from Corollary \ref{smooth:cor} that 
	\begin{equation}\label{gk:eq}
	\|g_m(A)J - J g_m(B)\|_{\GS}\le C (1+|m|)^{-\b} \|J\|^{1-\s} \||V|^\s\|_{\GS}, 
	\ \ m\in\mathbb Z,
	\end{equation}
	with a constant $C$ independent of $m$ and $g$.  
	Now use the $q$-triangle inequality \eqref{qtriangle:eq}: 
	\begin{align*}
		\|g(A) J - J g(B)\|_{\GS}^q\le &\ \sum_{m\in\mathbb Z} 
		\|g_m(A) J - J g_m(B)\|_{\GS}^q\\[0.2cm]
		\le &\ C\|J\|^{q(1-\s)} \||V|^\s\|_{\GS}^q
		\sum_{m\in\mathbb Z}(1+|m|)^{-q\b}.
	\end{align*}
	Since $q\b >1$, the above bound leads to \eqref{gest:eq}.
\end{proof}

\subsection{An important special case}
As explained in the Introduction, 
it is of particular interest for us to consider the case 
when the operator $J$ is an orthogonal projection. 

\begin{cond}\label{AP:cond}
Let $A$ be a self-adjoint operator on $\GH$, and let $P$ be an orthogonal projection 
such that $PD(A)\subset D(A)$, and $PA(I-P)$ extends to $\GH$ as a bounded operator.  
\end{cond}

The  condition $PD(A)\subset D(A)$ guarantees that $PAP$ is self-adjoint 
on the domain $PD(A)\oplus  (I-P)\GH$.

\begin{thm}\label{Szego1:thm} 
Suppose that $f$ 
satisfies Condition \ref{f:cond} with some $\g >0$,  $n\ge 2$ and $R > 0$. 
Let $\GS$ be a $q$-normed ideal where $(n-\s)^{-1} < q\le 1$ with some 
number $\s \in (0, 1]$, $\s < \g$.
Let $A, P$ be a self-adjoint operator and an orthogonal projection 
satisfying Condition \ref{AP:cond}.  
Suppose that 
$|PA(I-P)|^\s\in \GS$. Then 	
\begin{equation*}
 \|f(PAP)P - P f(A)\|_{\GS}
 \le C \1 f\1_n R^{\g-\s} \||PA(I-P)|^\s\|_{\GS}.
 \end{equation*}
 with a positive constant $C$ independent of the operators $A, P$, 
 function $f$ and parameter $R$.  
\end{thm}

\begin{proof} 
Denote 
\begin{equation*}
B_1 = PAP, B_2 = A,\ J = P.
\end{equation*}
Then
\begin{equation*}
f(PAP)P - P f(A) = f(B_1)J - J f(B_2).
\end{equation*}
Since $V = B_1J -JB_2 = - PA(I-P)$,  
Theorem \ref{quasinorm_est:thm} leads to the required estimate. 
\end{proof} 

We also state the following consequence of Corollary \ref{smooth:cor}:

\begin{cor}\label{smoothP:cor}
	Suppose that $g\in \plainC{n}_0(-\rho, \rho)$, with some $\rho >0$ and $n\ge 2$. 
	Let $\GS$ be a $q$-normed ideal with $(n-\s)^{-1} < q \le 1$, 
	where $\s\in (0, 1]$. 
	Let the operator $A$ and orthogonal projection $P$ be as in Theorem 
	\ref{Szego1:thm}. 
	Then 
	\begin{equation}\label{smooth:eq}
	\|g(PAP)P - P g(A)\|_{\GS}
	\le C  \max_{0\le k\le n} \bigl(\rho^k \| g^{(k)}\|_{\plainL\infty}\bigr)   
	\rho^{-\s}  \||PA(I-P)|^\s\|_{\GS},
	\end{equation}
	with a constant $C$ independent of the operator $A$ and projection $P$. 
\end{cor}

\section{Proof of Theorem \ref{quasinorm_est:thm}}

\subsection{A quasi-analytic extension}
In order to study functions of self-adjoint operators we use the formula known 
as the Helffer-Sj\"ostrand formula, see \cite{HelSjo}, \cite{EBD}. 
It requires the notion of 
a quasi-analytic extension of $f$. 
We use a somewhat more complicated definition than that 
in \cite{HelSjo} since we are working with non-smooth functions.  
For the sake of simplicity we concentrate on compactly supported functions, although all the 
definitions with appropriate modifications can be given for more general functions. 
Let 
\begin{equation*}
	\Pi = \Pi_+\cup \Pi_-,\  \Pi_{\pm} = \{z=(x, y): \pm y > 0\}.
\end{equation*} 
  
\begin{defn}\label{QA:defn} 
 Let 
 $f\in \plainC{0, \varkappa}_0(\R), 0 < \varkappa\le 1$,  
 and let 
$\tilde f\in \plainC{0, \varkappa}_0(\R^2)\cap\plainC1(\Pi)$  
be a function such that 
 \begin{enumerate}
\item
$\tilde f(x, 0) = f(x)$,  for all $x\in \R$, and  
\item
$|y|^{-1}\om\in \plainL1(\R^2)$, where 
\begin{equation*}
	\om(x, y) = \om(x, y; \tilde f) = \frac{\p}{\p \overline{z}} \tilde f(x, y) = \frac{1}{2}\biggl[
	\frac{\p}{\p x} \tilde f(x, y) + i\frac{\p}{\p y} \tilde f(x, y)
	\biggr].
\end{equation*}
\end{enumerate} 
Then $\tilde f$ is said to be a quasi-analytic extension of $f$. 
\end{defn}

 \begin{prop}\label{Hel_Sjo:prop}
 Let $A$ be a self-adjoint operator 
 on a Hilbert space $\GH$. Let $f:\R\to\R$
be a function as in Definition  \ref{QA:defn}, and let 
$\tilde f$ be its quasi-analytic extension. Then 
 \begin{equation}\label{Hel_Sjo:eq}
f(A) = \frac{1}{\pi}
\int\int  \frac{\p}{\p\overline z} \tilde f(x, y) R(x+iy; A) dx dy. 
 \end{equation}
 \end{prop}
 
 
 \begin{proof} 
 	It suffices to show that
 	\begin{equation*}
 f(t) = \frac{1}{\pi}
 \int\int  \frac{\p}{\p\overline z} \tilde f(x, y) (t- x-iy)^{-1} dx dy,\ \forall t\in\R. 
 	\end{equation*}
For a $\d>0$ split the plane into two regions:
\begin{equation*}
\CD_1 = \CD_1(\d) = \{z: |y|\ge \d\},
\ \ 
\CD_2 = \CD_2(\d) = \{z: |y|< \d\}.
\end{equation*} 	
First estimate the contribution 
from $\CD_2$:
\begin{equation*}
	\biggl| 
	\iint_{\CD_2} \omega(x, y) (t- x-iy)^{-1}   dx dy
		\biggr|
		\le \iint_{|y|<\d} |y|^{-1}|\omega(x, y)| dx dy.
\end{equation*} 	
 	By Definition \ref{QA:defn} the integral tends to zero 
 	as $\d\to 0$. Using the property 
 	\begin{equation*}
 		\frac{\p}{\p \overline{z}} \tilde f(x, y) (t - x-iy)^{-1}
 		= \frac{\p}{\p \overline{z}} \biggl(\tilde f(x, y) (t-x-iy)^{-1}\biggr),
 	\end{equation*}
 	we simplify the remaining integral:
 	\begin{equation*}
 	\iint_{\CD_1}  \omega(x, y) (t- x-iy)^{-1}   dx dy = 
 	\int_{\R}  F_{\d}(x; t)dx,
 	\end{equation*}
 	where
 	\begin{equation*}
F_{\d}(x; t) = \frac{1}{2\pi i} \bigl(\tilde f(x, \d)(t-x-i\d)^{-1} 
- \tilde f(x, -\d)(t-x+i\d)^{-1}\bigr).   		
 	\end{equation*}
 	Rewrite:
 	\begin{align*}
 2\pi i F_{\d}(x; t) = &
 f(x) \bigl((t-x-i\d)^{-1} - (t-x+i\d)^{-1}\bigr)\\[0.2cm] 
 + & \bigl(\tilde f(x, \d) - \tilde f(x, 0)\bigr)(t-x-i\d)^{-1}\\[0.2cm]
+ & \bigl(\tilde f(x, -\d) - \tilde f(x, 0)\bigr)(t-x+i\d)^{-1}.
\end{align*}
The last two terms converge to zero for all $x\not = t$ 
as $\d\to 0$. Moreover, since $f\in\plainC{0, \varkappa}$, we have 
\begin{equation*}
|\tilde f(x, \pm \d) - \tilde f(x, 0)|\le C\d^{\varkappa}, 
\end{equation*}
and hence the last two terms on the right-hand side do not exceed 
	$C |t-x|^{\varkappa-1}$.
Thus their integral over $x$ converges to zero as $\d \to 0$. Consequently,
\begin{equation*}
\int F_{\d}(x; t)dx
	= \frac{1}{2\pi i}\int f(x) \bigl((t-x-i\d)^{-1} - (t-x+i\d)^{-1}\bigr) dx 
	+ o(1),\ \d\to 0.
\end{equation*}
Since $f$ is continuous, the integral converges to $f(t)$, as claimed.
\end{proof}

Versions of the formula \eqref{Hel_Sjo:eq} have been known well before the paper 
\cite{HelSjo}. 
 In \cite{D1} E.M. Dyn'kin 
developed functional calculus for operators in Banach spaces, 
based on a formula in the spirit of \eqref{Hel_Sjo:eq}. 
Similar functional constructions can be found in L. H\"ormander's book 
\cite{Hor}, Section 3.1, so \eqref{Hel_Sjo:eq} must have been known to him 
earlier.  
In \cite{D2} E.M. Dyn'kin found a characterization 
of the classical Besov and Sobolev classes in terms of 
quasi-analytic extensions. These results were used in \cite{FP}.  
 
Let us describe a 
convenient quasi-analytic extension of the function 
$f$ satisfying Condition  
\ref{f:cond}. For convenience assume that $x_0 = 0$.  
For $b >0$ introduce the domain 
\begin{equation}\label{Fb:eq}
F_b = \{(x, y)\in\R^2: |y| < b|x|\}.
\end{equation}
By $U_b = U_b(x, y)$ we denote the characteristic function of $F_b$, i.e. 
\begin{equation}\label{ub:eq}
U_b(x, y) = 
\begin{cases}
1,\ |y|< b|x|,\\[0.2cm]
0,\  |y|\ge b|x|.
\end{cases}
\end{equation}
Let $\z\in\plainC\infty_0(\R)$ be a function 
such that 
\begin{equation}\label{zeta:eq}
\z(t) = 1 \ \ \textup{for}\ \  |t|\le 1/2, \ \ 
\textup{and}\ \  \z(t) = 0, |t|\ge 1.
\end{equation}

\begin{lem}\label{ab:lem}
Let $f$ satisfy Condition \ref{f:cond} 
with $x_0 = 0$ and some $n\ge 2$ and $R=1$. Then 
$f$ has a quasi-analytic extension 
$\tilde f\in\plainC{0, \varkappa}(\R^2)\cap \plainC1(\Pi)$, 
with the $\varkappa$ defined in \eqref{H:eq}, 
such that 
$\tilde f(x, y) = 0$ if $|y|>|x|$. Moreover, the derivative
\begin{equation*}
\om(x, y) = \frac{\p}{\p \overline{z}}\tilde f(x, y),
\end{equation*}
satisfies the bound
\begin{equation}\label{ab:eq}
|\om(x, y)|
\le C_n \1 f\1_n |x|^{\g-n} |y|^{n-1} U_1(x, y) \chi_1(x),
\end{equation}
for $x\not = 0, y\in\R$. 
The constant $C_n$ does not depend on $f$. 
\end{lem}

\begin{proof} 	
	We use a slight modification of the 
	``standard" construction of a quasi-analytic extension 
which can be found e.g. in \cite{Hor}, Section 3.1, 
or \cite{EBD}, Ch. 2. 
Let the function $\z$ 
be defined as in \eqref{zeta:eq}. Without loss of generality assume $\1 f\1_n = 1$.   
Define for all $x\not = 0$ and $y\in\R$: 
\begin{equation*} 
\tilde f(x, y) = \biggl[\sum_{l=0}^{n-1} 
f^{(l)}(x) \frac{(iy)^l}{l!}\biggr]\s(x, y),\ 
\s(x, y) = 
 \z\biggl(\frac{y}{ x }\biggr).
\end{equation*} 
For $x = 0$ we set $\tilde f(0, y) = 0, y\in \R$. 
Clearly, $\tilde f(x, 0) = f(x)$, for all $x\in \R$. 
Moreover,  $\tilde f$ is trivially continuous for all $x\not =0$. 
At $x = 0$ it is continuous because of the bound 
$|\tilde f(x, y)|\le C|x|^\g$, for all 
$x\not = 0, y\in\R$, which follows from 
\eqref{fbound:eq}. Furthermore, $\tilde f\in \plainC{1}(\mathbb C\setminus\{0\})$ and  
one checks directly for $x\not = 0$ that 
\begin{equation*}
\frac{\p \tilde f}{\p x}(x, y) = 
\biggl[\sum_{l=0}^{n-1} f^{(l+1)}(x) \frac{(iy)^l}{l!}\biggr] \s(x, y) 
+ \biggl[\sum_{l=0}^{n-1} f^{(l)}(x) \frac{(iy)^l}{l!}\biggr] \frac{\p\s}{\p x}(x, y),
\end{equation*}
and 
\begin{equation*}
\frac{\p \tilde f}{\p y}(x, y) = 
i\biggl[\sum_{l=1}^{n-1} f^{(l)}(x) \frac{(iy)^{l-1}}{(l-1)!}\biggr] \s(x, y) 
+ \biggl[\sum_{l=0}^{n-1} f^{(l)}(x) \frac{(iy)^l}{l!}\biggr] \frac{\p\s}{\p y}(x, y).
\end{equation*}
 Thus 
\begin{align}
\om(x, y) = \frac{\p}{\p\overline z} \tilde f(x, y)
= &\ \frac{1}{2}\biggl( \frac{\p \tilde f}{\p  x}  
+ i\frac{\p\tilde f}{\p y}\biggr)(x, y) \notag\\[0.2cm]
= &\ \frac{1}{2} f^{(n)}(x) \frac{(iy)^{n-1}}{(n-1)!}\s(x, y)
+ \biggl[\sum_{l=0}^{n-1} f^{(l)}(x) \frac{(iy)^l}{l!}\biggr]  
\frac{\p\s}{\p \overline z}(x, y),
\label{deriv:eq}
\end{align}
for all $x\not = 0$.
Since we have assumed $\1 f\1_n = 1$, 
\eqref{fbound:eq} implies that for $x\not = 0$ we have 
\begin{equation}\label{ab3:eq}
|\tilde f_y(x, y)|
\le C_n  
\biggl[ \sum_{l=1}^{n-1}|x|^{\g-l} |y|^{l-1} |\z(yx^{-1})| 
+ |x|^{-1}\biggl[\sum_{l=0}^{n-1} |x|^{\g-l} |y|^l\biggr]|\z'(yx^{-1})|\biggr]
\chi_1(x). 
\end{equation}
Recall that $|y| < |x|\le 1$ on the support of $\zeta(yx^{-1})$, so that
\begin{equation*}  
|\tilde f_y(x, y)|\le C|x|^{\gamma-1} U_1(x, y)\chi_1(x), x\not = 0.
\end{equation*}
Using the bound $|y|<|x|\le1$ again, we deduce  
that   
\begin{equation*}
|\tilde f_y(x, y)|\le C|y|^{\varkappa-1} U_1(x, y)\chi_1(x), y\not = 0.	
\end{equation*}  
where $\varkappa = \min\{\gamma, 1\}$. 
One easily proves the same bounds for the derivative 
$\tilde f_x$. These bounds imply that $\tilde f\in \plainC{0, \varkappa}(\R^2)$.  

In order to establish \eqref{ab:eq} we use \eqref{deriv:eq}, so that 
for all $x\not = 0$ we have 
\begin{equation}\label{ab1:eq}
|\om(x, y)|
\le C_n  
\biggl[ |x|^{\g-n} |y|^{n-1} |\z(yx^{-1})| 
+ |x|^{-1}\biggl[\sum_{l=0}^{n-1} |x|^{\g-l} |y|^l\biggr]|\z'(yx^{-1})|\biggr]
\chi_1(x). 
\end{equation}
The first term on the right-hand side already satisfies \eqref{ab:eq}.  
Using the formula 
\begin{equation*} 
\z'(yx^{-1}) = \z'(yx^{-1}) U_1(x, y) \bigl(1- U_{\frac{1}{2}}(x, y)\bigr),
\end{equation*} 
we conclude that in the second term on the right-hand side 
of \eqref{ab1:eq} we have $|x|/2\le |y|\le |x|$. Hence this term satisfies 
\eqref{ab:eq} as well. 

Since $\g>0$, the bound \eqref{ab:eq} ensures that 
$|y|^{-1}\om\in\plainL1(\R^2)$. 
This completes the proof. 
\end{proof}

\subsection{Proof of Theorem \ref{quasinorm_est:thm}}
Without loss of generality we may assume that $x_0 = 0$, 
and that $\1 f\1_n = 1$. 
Also, in view of Remark \ref{fest:rem}, it suffices to obtain \eqref{fest:eq} 
for $R=1$ only. 

Let $\tilde f$ be the quasi-analytic extension  
constructed in Lemma \ref{ab:lem}. 
By the formula \eqref{Hel_Sjo:eq} and 
resolvent identity \eqref{res:eq} 
we can write 
\begin{align*}
T = f(A)J- Jf(B)
= &\ - \frac{1}{\pi}\iint \om(x, y) R(z; A) V R(z; B) dx dy, z = x+iy.
\end{align*}
Let $\rho(x, y) = 8^{-1}|y|$, and let $\{\CD_j\}, j=1, 2, \dots,$ 
be a family of open discs with finite intersection property 
centred at some points $z_j = (x_j, y_j)\in \Pi$, 
of radius $\rho(x_j, y_j)$ such that 
\begin{equation*}
\cup_j \CD_j = \Pi.
\end{equation*}
Let $\phi_j\in \plainC\infty_0(\Pi)$ be an associated partition of unity 
such that 
\begin{equation*}
|\p_x^l \p_y^k \phi_j(x, y)|\le C_{l, k} \rho(x, y)^{-l-k}.  
\end{equation*} 
By the ``finite intersection property" we mean that 
the number of discs having non-empty common intersection is uniformly bounded. 
The existence of such a covering and such a partition of unity follows from 
\cite{Hor}, Theorem 1.4.10. 
Estimate the quasi-norm of
\begin{equation*}
T_j =  \iint  \phi_{j}(x, y)\om(x, y) R(z; A) V R(z; B) dx dy.
\end{equation*}
For $z\in \CD_j$ 
expand $R(z; A)$ 
and 
$R(z; B)$ in the uniformly norm-convergent series 
\begin{equation*}
R(z; K) = \sum_{k=0}^\infty (z-z_j)^k R(z_j; K)^{k+1},\ 
\end{equation*}
where $K=A$ or $B$. 
The uniformity of convergence is guaranteed by the bound  
$|z-~z_j| \|R( z_j; K)\|\le 1/8$. 
Denote $\bk = (k_1, k_2),\ k_1\ge 0, k_2\ge 0$. 
Therefore 
%
we can now expand $T_j$ in the norm-convergent series:
\begin{align*}
T_j = \sum_{\bk} T_{j\bk},\ 
T_{j\bk} = &\ \biggl[\iint  \om(x, y)\phi_{j}(x, y) (z-z_j)^{k_1+k_2}  
dx dy \biggr] R(z_j; A)^{k_1} W_j R(z_j; B)^{k_2},\\[0.2cm]
W_j = &\ R(z_j; A) V R(z_j; B).
\end{align*} 
By Lemma \ref{sandwich:lem}, 
\begin{equation*}
\|W_j\|_{\GS}\le 2^{1-\s}\|J\|^{1-\s} \||V|^\s\|_{\GS}|y_j|^{-1-\s}.
\end{equation*}
Estimate the $\GS$-quasi-norm of each $T_{j\bk}$. It follows from 
\eqref{ab:eq} that 
\begin{equation*}
\|T_{j\bk}\|_{\GS}\le 
C  |y_j|^{-1-\s} 8^{-k_1+k_2}\|J\|^{1-\s}\ \||V|^\s\|_{\GS} 
\underset{\CD_j}\iint |x|^{\g-n}|y|^{n-1} U_1(x, y)\chi_1(x)dx dy, 
\end{equation*} 
A straightforward calculation shows that if $\CD_j\cap F_1\not = \varnothing$ then 
$(x_j, y_j)\in F_2$ and $\CD_j\subset F_4$, see 
\eqref{Fb:eq} for the definition of $F_b, b >0$.
Thus for all $(x, y)\in\CD_j$ we have 
\begin{equation*}
|x-x_j|< 8^{-1} |y_j| < 4^{-1}|x_j|,\quad\  |y-y_j|< 8^{-1} |y_j|,
\end{equation*}
so
\begin{equation*}
\frac{3}{4}|x_j| < |x| < \frac{5}{4}|x_j|,\quad\  
\frac{7}{8}|y_j| < |y| < \frac{9}{8}|y_j|.
\end{equation*}
Consequently,
\begin{equation*}
\|T_{j\bk}\|_{\GS}\le 
C   8^{-k_1+k_2}\|J\|^{1-\s}\  \| |V|^{\s} \|_{\GS}  |x_j|^{\g-n}|y_j|^{n-\s} 
U_2(x_j, y_j) \chi_{2}(x_j).  
\end{equation*}
Since $\CD_j\subset F_4$, by the $q$-triangle inequality, we have 
\begin{align*}
\|T_j\|_{\GS}^q\le &\ \sum_{\bk} \|T_{j\bk}\|_{\GS}^q\\[0.2cm]
\le &\ C  \|J\|^{q(1-\s)}\  \||V|^\s\|_{\GS}^q \  
|x_j|^{q(\g - n)}|y_j|^{q(n-\s)} U_2(x_j, y_j) \chi_{2}(x_j)
\sum_{\bk} 8^{-(k_1+k_2)q}\\[0.2cm]
\le &\ C   \|J\|^{q(1-\s)}\  \||V|^\s\|_{\GS}^q \  
\underset{\CD_j}\iint |x|^{q(\g-n)}|y|^{q(n-\s)-2} U_4(x, y) \chi_{4}(x)dx dy.
\end{align*}
Use the $q$-triangle inequality again to sum over $j$:
\begin{align*}
\|T\|_{\GS}^q\le &\ \sum_{j}\|T_j\|_{\GS}^q\\[0.2cm]
\le &\ 
C  \|J\|^{q(1-\s)}\||V|^\s\|_{\GS}^q \ 
\sum_j \underset{\CD_j}\iint |x|^{q(\g-n)}|y|^{q(n-\s)-2} U_4(x, y) 
\chi_{4}(x)dx dy\\[0.2cm]
\le &\  
C  \|J\|^{q(1-\s)}\||V|^\s\|_{\GS}^q 
\iint |x|^{q(\g - n)}|y|^{q(n - \s) - 2} U_4(x, y) \chi_{4}(x)dxdy,
\end{align*}
where we have used the finite intersection property. 
By assumption, we have $q(n-\s)-2 > -1$, and hence 
the integral on the right-hand side is bounded by
\begin{equation*}
C\int_{|x|<4} |x|^{q(\g-n)} 
\biggl[\int_{|y|\le 4|x|} |y|^{q(n-\s) - 2}dy\biggr] dx
\le \tilde C \int_0^{4} x^{q(\g-\s)-1}dx \le  C. 
\end{equation*}
This proves \eqref{fest:eq} for $R = 1$. 
As explained in Remark \ref{fest:rem}, this immediately leads to 
\eqref{fest:eq} for general $R>0$. 
\qed

\section{Trace asymptotics for multidimensional Wiener-Hopf operators with discontinuous symbols} 

\subsection{Definitions}
Now we derive from the theorems established above estimates for 
some Wiener-Hopf operators on $\plainL2(\R^d)$. 
Let $\L, \Om\subset \R^d, d\ge 2,$ be two domains, 
and let $\chi_\L$, $\chi_\Om$ be their characteristic functions.   
For a bounded complex-valued function $a=a(\bx, \bxi)$, called \textit{symbol}, 
define the pseudo-differential operator 
\begin{equation*}
\bigl(\op_\a(a) u\bigr)(\bx)
= \frac{\a^{d}}{(2\pi)^{\frac{d}{2}}}
\iint e^{i\a\bxi\cdot(\bx-\by)} a(\bx, \bxi) u(\by) d\by d\bxi,  u\in \plainS(\R^d).
\end{equation*}
It is a standard fact that 
under the condition $a\in \plainW{d+1}{\infty}(\R^d\times\R^d)$  
the norm of the operator 
$\op_\a(a)$ is bounded uniformly in $\a\ge 1$, see e.g. 
\cite{Sob}, Lemma 3.9.

Under  Wiener-Hopf operators with discontinuous symbols 
here we understand operators of the form 
\begin{equation}\label{WH:eq}
S_\a(a) = S_\a(a; \L, \Om) 
= \chi_\L P_{\Om, \a} 
\re \op_\a(a) P_{\Om, \a}\chi_\L,\ \  \ P_{\Om, \a} = \op_\a(\chi_\Om).
\end{equation}
The function $a$ is assumed to be smooth, and it is the presence of 
the projection $P_{\Om, \a}$ that suggests the term 
``discontinuous symbol". In Sect. \ref{more:sect} we consider somewhat more 
general discontinuous symbols. 

We impose the following conditions on the domains $\L$ and $\Om$.  

\begin{cond}\label{domain:cond}
The domains $\L, \Om\subset \R^d, d\ge 2$, are both Lipschitz domains; 
$\Om$ is bounded, 
and either $\L$ or $\R^d\setminus\L$ is bounded. 
\end{cond}

\begin{cond}\label{smooth:cond}
The domain $\L$ is piece-wise $\plainC1$, and 
$\Om$ is piece-wise $\plainC3$.
\end{cond}

By the Lipschitz domain we understand a domain which locally looks like a 
set of points 
above the graph of a suitable Lipschitz function. 
Precise definitions of this property as well as of 
piece-wise smoothness are given in \cite{Sob2}, 
Definition 2.1. 

We are interested in the large $\a$ asymptotics of the trace of the operator 
\begin{equation}\label{Dal:eq}
D_\a(a, \L, \Om; g)
= \chi_\L g\bigl(S_\a(a, \L, \Om)\bigr)\chi_\L 
- \chi_\L g\bigl( S_\a(a, \R^d, \Om)\bigr)\chi_\L
\end{equation}
with a function $g:\R\to\mathbb C$ 
which is smooth except for finitely many points.

Let us define the asymptotic coefficients entering the main asymptotic formulas. 
For a symbol $b=b(\bx, \bxi)$ let 
\begin{equation}\label{w0:eq}
\GW_0(b) = \GW_0(b; \L, \Om)
= \frac{1}{(2\pi)^d}\int_\Om\int_{\L} b(\bx, \bxi) d\bx d\bxi.
\end{equation}
For any $(d-1)$-dimensional Lipschitz surfaces $L, P$ denote
\begin{equation}\label{w1:eq}
\GW_1(b) = \GW_1(b ; L, P) = \frac{1}{(2\pi)^{d-1}}\int_{L} \int_{P} b(\bx, \bxi)
| \bn_{L}(\bx)\cdot\bn_{P}(\bxi) |  dS_{\bxi} dS_{\bx},
\end{equation}
where $\bn_{L}(\bx)$ and $\bn_P(\bxi)$ denote the exterior unit normals to
$L$ and $P$ defined for a.e. $\bx$ and $\bxi$ respectively. 
For any function $g\in \plainC{0, \varkappa}(\mathbb C), \varkappa>0$, 
and any number $s\in\mathbb C$, 
we also define
\begin{equation}\label{GA:eq}
\GA(g; s) = \frac{1}{(2\pi)^2}\int_0^1 
\frac{g(st) - (1-t) g(0)- t g(s)}{t(1-t)} dt.
\end{equation}
The next result is found in \cite{Sob2}, Theorem 2.5:

\begin{prop}\label{Widom:prop}
Let $\L, \Om\subset\R^d, d\ge 2,$  be two domains satisfying Conditions 
\ref{domain:cond} and \ref{smooth:cond}. 
Let $a\in \plainW{d+2}{\infty}(\R^d\times\R^d)$ be a complex-valued function. 
Let $g_p(t) = t^p$ with some $p = 0, 1, \dots$. 
Then 
\begin{equation}\label{Widom:eq}
\lim_{\a\to\infty}
\frac{1}{\a^{d-1}\log\a}\tr D_\a(a, \L, \Om; g_p )
= \GW_1(\GA(g_p; \re a), \p\L, \p\Om).
\end{equation}
\end{prop}

Observe that for the polynomials $g_0(t)\equiv 1$ and $g_1(t) = t$ 
both sides of 
the above formula equal zero, so the asymptotics hold trivially. 

The main focus of of \cite{Sob2} was on non-smooth 
domains $\L$ and $\Om$. As a result 
the above proposition was formally proved in \cite{Sob2} 
for the case $d\ge 2$ only, although a similar, somewhat simplified 
argument should give \eqref{Widom:eq} for the case $d=1$ as well, with an appropriately 
modified definition of the coefficient $\GW_1$, see e.g. 
\cite{Widom_82}. However we do not pursue this objective in the 
current paper.  

Our aim here is to extend Proposition 
\ref{Widom:prop} to functions $g$ that 
have just $\plainC2$ local smoothness, and may lose differentiability 
at finitely many points.

\begin{thm}\label{Szego3:thm} Let $d\ge 2$, and 
let the domains $\L, \Om$ 
be two domains satisfying Conditions 
\ref{domain:cond} and \ref{smooth:cond}. 
Let the symbol 
$a = a(\bx, \bxi)$ be a globally 
bounded $\plainC\infty$-function.  
Let $X = \{z_1, z_2, \dots, z_N\}\subset \R$, $N <\infty$, be a collection of 
points on the real line. 
Suppose that $g\in\plainC2(\R\setminus X)$ is a function such that 
in a neighbourhood of each point $z\in X$ 
it
satisfies the bound
\begin{equation*}
|g^{(k)}(t)|\le C_k |t - z|^{\g-k}, \ k = 0, 1, 2, 
\end{equation*} 
with some $\g>0$. Then 
\begin{equation}\label{Widom_new:eq}
\lim_{\a\to\infty}\frac{1}{\a^{d-1}\log\a} \tr D_\a(a, \L, \Om; g) 
= \GW_1\bigl(
\GA(g; \re a); \p\L, \p\Om
\bigr).
\end{equation}
\end{thm}

Let us make some comments on Theorem \ref{Szego3:thm}.

\begin{rem} \label{expl:rem}

\begin{enumerate}  
\item
The assumption $a\in \plainC\infty$ is made for simplicity. In fact, 
some finite smoothness, depending on the value of the parameter $\gamma$, would 
suffice, but we have chosen to avoid ensuing technicalities. 	 
	
\item\label{bounded:item}
Suppose that $\L$ is bounded 
and that $g\in \plainC\infty(\R)$ is a function such that $g(0) = 0$. Then 
both operators on the right-hand side of \eqref{Dal:eq} are trace class, and 
formula \eqref{Widom_new:eq} is just an indirect way to write the asymptotics 
\begin{align}\label{Widom1:eq}
	\tr g(S_\a) = \a^d &\ \GW_0(g(\re a); \L, \Om)\notag\\[0.2cm]
	&\ + \a^{d-1} \log\a \ 
	\GW_1(\GA(g; \re a); \p\L, \p\Om) + o(\a^{d-1}\log\a),
\end{align}
as $\a\to\infty$, established in \cite{Sob2}, Theorem 2.3. Indeed, 
one can show that 
\begin{equation*}
\tr \chi_\L g(S_\a(a; \R^d, \Om))\chi_\L = \a^d \GW_0(g(\re a); \L, \Om) + O(\a^{d-1}), \a\to\infty,   
\end{equation*}
see e.g. \cite{Sob}, Section 12.3, where a similar calculation was done. 
Substituting this in \eqref{Widom_new:eq} one obtains \eqref{Widom1:eq}. 

If the symbol $a$ depends only on the 
variable $\bxi$, i.e. $a = a(\bxi)$, then the 
reduction of \eqref{Widom_new:eq} to \eqref{Widom1:eq} for bounded $\L$ 
is more straightforward. 
Indeed, in this case the 
second operator on the right-hand 
side of \eqref{Dal:eq} is given by 
\begin{equation}\label{second:eq}
\chi_\L \op_\a \bigl(g(\re a\chi_\Om)\bigr) \chi_\L.
\end{equation}
This operator is clearly trace class for any continuous $g$ such that $g(0) = 0$.  
Integrating 
its kernel along the diagonal, 
one easily finds the exact value of its trace: 
$\a^d\GW_0(g(\re a); \L, \Om)$.


\item 
Suppose that $a = a(\bxi)$, and that $g=0$ on the range of the function $\re a$. 
Then the operator \eqref{second:eq} equals zero, 
so that Theorem 4.4 implies that $g(S_\a(a,\L, \Om))$ is trace class and its trace 
satisfies the asymptotic formula
\begin{equation}\label{Widom_new_new:eq}
	\lim_{\a\to\infty}\frac{1}{\a^{d-1}\log\a} \tr g\bigl(S_\a(a, \L, \Om)  \bigr)
	= \GW_1\bigl( \GA(g; \re a); \p\L, \p\Om \bigr),
\end{equation}
It is interesting to point out that 
this formula holds for both bounded or unbounded domains $\L$. 

Another proof of formula 
\eqref{Widom_new_new:eq} in 
the special case $a \equiv 1$ and $g(0) = g(1) = 0$ was given in \cite{LeSpSo}. 
It was motivated by the study of the 
\textit{entanglement entropy} for free  Fermions 
at zero temperature, see also \cite{GiKl} and \cite{HLS}. 
Mathematically speaking, the entropy is found as trace of the operator 
$\eta_\b(S_\a), \b >0$,  
with the operator $S_\a = S_\a(1, \L, \Om)$ for bounded $\L$ and $\Om$, 
and with  the function $\eta_\b $ defined by 
\begin{equation}\label{eta:eq}
\eta_\b(t) =
\begin{cases}
\dfrac{1}{1-\b} \log(t^\b + (1-t)^\b),\   \b>0, \b\not =1 \ (\textup{R\'enyi entropy}),\\[0.4cm]
-t\log t - (1-t)\log (1-t),\ \b = 1	(\textup{von Neumann entropy}),
\end{cases}
\end{equation} 
if $t \in [0, 1]$, end extended 
by $0$ to the rest of the real line. 
Clearly, for $\b \not = 1$ the function  
$\eta_\b$ satisfies the conditions of Theorem \ref{Szego3:thm} 
with $\g = \b$, and for $\b  = 1$ -- with arbitrary $ \g < 1$.

 \end{enumerate}

\end{rem}

Before proving Theorem \ref{Szego3:thm} we list some useful facts. 

\begin{lem} 
\begin{enumerate}
\item
If $g\in \plainW{1}{\infty}(\R)$,
then 
\begin{equation}\label{diff:eq}
|\GA(g; s)|\le \frac{1}{\pi^2}|s| \|g'\|_{\plainL\infty}.
\end{equation}
\item
Suppose that $g$ satisfies \eqref{fbound:eq} with some $0 < R\le 1$ . 
Then 
\begin{equation}\label{sing:eq}
|\GA(g; s)|\le  
C |s|^{\frac{\varkappa}{2}}R^{\frac{\g}{2}} \1 g\1_1,\ \varkappa = \min\{1, \g\}.
\end{equation}
\end{enumerate}
\end{lem}

\begin{proof} 
Using the formula
$t^{-1}(1-t)^{-1} = t^{-1}+(1-t)^{-1}$, 
we rewrite $\GA$ in the form
\begin{equation}\label{newa:eq}
(2\pi)^2\GA(g; s) = \int_0^1\frac{g(st)-g(0)}{t} dt 
+ \int_0^1\frac{g(st)-g(s)}{1-t} dt.
\end{equation}
Now \eqref{diff:eq} follows from the bounds
\begin{equation*}
|g(st)-g(0)|\le \|g'\|_{\plainL\infty} |s||t|,\ \ 
 |g(st)-g(s)|\le \|g'\|_{\plainL\infty} |s||t-1|.
\end{equation*}

Proof of \eqref{sing:eq}. 
Assume without loss of generality that $\1 g\1_1 = 1$.  
By \eqref{H:eq} 
the first integral in \eqref{newa:eq} 
is estimated by
\begin{equation*}
\bigl(2\|g\|_{\plainL\infty}\bigr)^{\frac{1}{2}}
	\int_0^1 \frac{|g(st)-g(0)|^{\frac{1}{2}}}{t}dt
	\le 
	C R^{\frac{\g}{2}}|s|^{\frac{\varkappa}{2}}
\int_0^1 t^{\frac{\varkappa}{2}-1} dt
\le \tilde C R^{\frac{\g}{2}}|s|^{\frac{\varkappa}{2}}. 
\end{equation*}
Similarly, the second integral is bounded by 
\begin{equation*}
	\bigl(2\|g\|_{\plainL\infty}\bigr)^{\frac{1}{2}}
	\int_0^1 \frac{|g(st)-g(s)|^{\frac{1}{2}}}{1-t}dt
	\le 
	C R^{\frac{\g}{2}}|s|^{\frac{\varkappa}{2}}\int_0^1 (1-t)^{\frac{\varkappa}{2}-1} dt
	\le \tilde C R^{\frac{\g}{2}}|s|^{\frac{\varkappa}{2}}. 
\end{equation*}
This proves \eqref{sing:eq}.
\end{proof}

A crucial ingredient in the proof of Theorem 
\ref{Szego3:thm} 
is the following lemma. 

\begin{lem}\label{cross:lem}
Let $\L$ and $\Om$ satisfy Condition 
\ref{domain:cond}, and let the symbol $a$ be as in Theorem \ref{Szego3:thm}. 
Then for any $q\in (0, 1]$, and all $\a\ge 2$ we have 
\begin{equation*}
\|\chi_\L P_{\Omega, \a}\op_\a(a) P_{\Om, \a}(I-\chi_\L)\|_{\GS_q}^q\le C_q\a^{d-1}\log\a,
\end{equation*}
with a constant $C_q$ independent of $\a$. 
\end{lem}

The above bound can be derived from \cite{Sob1} Theorem 4.6 in the same
way as \cite{Sob1} Corollary 4.7.

\subsection{ 
Proof of Theorem \ref{Szego3:thm}} 
Throughout the proof we denote for brevity
$D_\a(g) = D_\a(a, \L, \Om; g)$, and 
$\GW_1(g) = \GW_1(\GA(g; \re a); \p\L, \p\Om)$. 
Recall that the norm of the operator 
$\op_\a(a)$ is bounded uniformly in $\a\ge 1$, so  
without loss of generality we may assume that  
~$\|\op_\a(a)\|\le ~1/2$ and that $g$ 
is supported on the interval $[-1, 1]$, and it is real-valued. 

The proof splits into two parts. 

\underline{Step 1: proof of formula \eqref{Widom_new:eq}  
for $g\in \plainC2(\R)$.} 
By the Weierstrass Theorem, for any $\varepsilon>~0$ 
one can find a real polynomial 
$g_\varepsilon$ such that the function 
$f_\varepsilon = g-g_\varepsilon$ satisfies the bound 
\begin{equation}\label{feps:eq}
\max_{0\le k\le 2}\max_{ |t|\le 1}|f_\varepsilon^{(k)}(t)|< \varepsilon.
\end{equation}
Now we use Corollary \ref{smoothP:cor} 
with $\GS = \GS_1$, $n=2$, $R=1$, arbitrary $\s \in (0, 1)$, 
and the operators
\begin{equation}\label{subst:eq}
A = P_{\Om, \a} \re \op_\a(a) P_{\Om, \a},\  P = \chi_\L.
\end{equation}
Corollary \ref{smoothP:cor}, Lemma \ref{cross:lem} and 
bound \eqref{feps:eq} give that 
\begin{align*}
\|D_\a(f_\varepsilon)\|_{\GS_1}
\le &\ C \varepsilon  
\| \chi_\L P_{\Om, \a}
\op_\a(a) P_{\Om, \a}(I-\chi_\L)\|_{\GS_\s}^\s\\[0.2cm]
\le &\ C \varepsilon \a^{d-1} \log \a, \ \a\ge 2,
\end{align*}
and as a consequence,
\begin{equation*}
\tr D_\a(g)\le \tr D_\a(g_\varepsilon) + \|D_\a(f_\varepsilon)\|_{\GS_1}
\le \tr D_\a(g_\varepsilon) + 
C\varepsilon\a^{d-1}\log\a.
\end{equation*}
Now, using Proposition 
\ref{Widom:prop} for the polynomial $g_\varepsilon$ we get
\begin{equation*}
\limsup_{\a\to\infty }\frac{1}{\a^{d-1}\log\a}
D_\a(g) 
\le \GW_1(g_\varepsilon)+ C\varepsilon.
\end{equation*}
Due to \eqref{diff:eq} and \eqref{feps:eq}, 
the asymptotic coefficient $\GW_1(f_\varepsilon)$ tends 
to zero as $\varepsilon\to 0$, so that 
\begin{equation*}
	\GW_1(g_\varepsilon) = \GW_1(g) - \GW_1(f_\varepsilon) \to \GW_1(g),\ 
	\varepsilon\to 0. 
\end{equation*}
This implies that 
\begin{equation*}
\limsup_{\a\to\infty }\frac{1}{\a^{d-1}\log\a}
D_\a(g) \le \GW_1(g).
\end{equation*}
In the same way one obtains the appropriate lower 
bound for the $\liminf$. This completes the proof of \eqref{Widom_new:eq} 
for $g\in \plainC2(\R)$. 

\underline{Step 2. Completion of the proof.} 
Let $g$ be a function as specified in the theorem. As before we assume that 
$g$ is real-valued and that it is supported on the interval $[-1, 1]$. 
By choosing an appropriate partition of unity we may assume that the set 
$X$ consists of one point only, which we denote by $z$. 
 
Let $\z\in\plainC\infty_0(\R)$ be a real-valued function, satisfying \eqref{zeta:eq}.
Represent $g= g_R^{(1)}+ g_R^{(2)}, 0<R\le 1$, where 
$g_R^{(1)}(t) = g(t) \z\bigl((t-z)R^{-1}\bigr)$,\ 
$g_R^{(2)}(t) = g(t)  - g_R^{(1)}(t)$. 
It is clear that  $g_R^{(2)}\in \plainC2(\R)$, 
so one can use the formula \eqref{Widom_new:eq} established in  
Part 1 of the proof:
\begin{equation}\label{g2R:eq}
\lim_{\a\to\infty }\frac{1}{\a^{d-1}\log\a}
D_\a(g_R^{(2)}) 
=  \GW_1(g_R^{(2)}).
\end{equation}
For $g_R^{(1)}$ we use Theorem 
\ref{Szego1:thm} with $\GS = \GS_1$, $n=2$, an arbitrary $\s\in (0, 1), \s <\g$,  
and with the operators $A, P$ defined in \eqref{subst:eq}. Noticing that 
$\1 g_R^{(1)}\1_2\le C\1 g\1_2$,  we get from Theorem \ref{Szego1:thm} 
and Lemma \ref{cross:lem} that 
\begin{equation*}
\|D_\a(g_R^{(1)})\|_{\GS_1}\le C_\s R^{\g-\s} 
\|(I-\chi_\L)A \chi_\L\|_{\GS_\s}^\s\le C_\s R^{\g-\s} \a^{d-1}\log\a,\ C_\s = C_\s(g), 
\end{equation*}
for all $\a \ge 2$. Therefore
\begin{equation*}
\tr D_\a(g)\le \tr D_\a(g_R^{(2)}) + C_\s R^{\g-\s} \a^{d-1}\log\a. 
\end{equation*}
Using \eqref{g2R:eq} we obtain the bound 
\begin{equation*}
\limsup_{\a\to\infty }\frac{1}{\a^{d-1}\log\a}
D_\a(g) 
\le \GW_1(g_R^{(2)})+ CR^{\g-\s}.
\end{equation*}
Due to \eqref{sing:eq} 
the asymptotic coefficient $\GW_1(g^{(1)}_R)$ converges to zero as $R\to 0$.  Thus
\begin{equation*}
\GW_1(g^{(2)}_R) = \GW_1(g) - \GW_1(g^{(1)}_R)\to \GW_1(g),\ R\to 0.
\end{equation*}
This implies that 
\begin{equation*}
\limsup_{\a\to\infty }\frac{1}{\a^{d-1}\log\a}
D_\a(g) \le \GW_1(g).
\end{equation*}
In the same way one obtains the appropriate lower 
bound for the $\liminf$. This completes the proof.
\qed

\section{More on symbols with jump discontinuities}\label{more:sect}

In this section we give a variant of Theorem \ref{Szego3:thm} 
%
%
for more general discontinuous symbols: instead of the symbols 
$a(\bx, \bxi)\chi_\Om(\bxi)$ we study 
$a(\bx, \bxi)\chi_\Om(\bxi) + a_1(\bx, \bxi)\chi_{\Om_1}(\bxi)$,\ 
$\Om_1 = \Rd\setminus\Om$, 
i.e  symbols allowed to have jump discontinuities on the surface $\p\Om$.

%
Along with the operator \eqref{WH:eq} introduce the notation for its 
non-symmetric variant: 
\begin{equation*} 
T_\a(a) = T_\a(a; \L, \Om) 
= \chi_\L P_{\Om, \a} 
\op_\a(a) P_{\Om, \a}\chi_\L,
\end{equation*}
so that $S_\a(a) = \re T_\a(a)$. 
Let $\Om_1 = \Rd\setminus\Om$, and let $a, a_1$ be two smooth symbols.  
Define 
\begin{align*}
	V_\a(a, a_1) = &\ V_\a(a, a_1; \L, \Om) = 
	T_\a(a; \L, \Om) + T_\a(a_1; \L, \Om_1), \\[0.2cm]
	H_\a(a, a_1) = &\ H_\a(a, a_1; \L, \Om)
	= S_\a(a; \L, \Om) + S_\a(a_1; \L, \Om_1).
\end{align*}
%
%
Both symbols $a$ and $a_1$ are assumed to 
have compact supports in the variable $\bxi$.

To state the result we need to define instead of \eqref{GA:eq} the coefficient 
\begin{equation}\label{GD:eq}
\GD(g; s, s_1) = \frac{1}{(2\pi)^2}\int_0^1 
\frac{g(st + s_1(1-t)) - t g(s)- (1-t) g(s_1)}{t(1-t)} dt, 
\end{equation}
$s, s_1\in\mathbb C$. 

\begin{thm}\label{Szego4:thm} 
Let the domains $\L, \Om$ and  function $g$ be in Theorem 
\ref{Szego3:thm}. 
Suppose that the symbols $a, a_1\in \plainC\infty(\Rd\times\Rd)$ 
 are globally bounded functions compactly supported in the variable $\bxi$. 
Then 
\begin{align}\label{jump:eq}
	\lim_{\a\to\infty}\frac{1}{\a^{d-1}\log\a}
	\tr\bigl[
	&\ \chi_\L g\bigl(H_\a(a, a_1; \L, \Om)\bigr) \chi_\L
	- \chi_\L g\bigl( H_\a(a, a_1; \Rd, \Om)\bigr)\chi_\L
	\bigr]\notag\\[0.2cm]
	&\ \quad\quad  =  \GW_1\bigl(
	\GD(g; \re a, \re a_1); \p\L, \p\Om \bigr).
\end{align}
\end{thm}

It is appropriate to make a comment in the spirit of Remark \ref{expl:rem} (\ref{bounded:item}):

	If the domain $\L$ is bounded and $g\in \plainC\infty(\R)$ is such that 
	$g(0) = 0$, then both operators on the left-hand side of 
	\eqref{jump:eq} are trace class, and formula \eqref{jump:eq} 
	is just another way to write the asymptotics 
	\begin{align*}
\tr g\bigl(H_\a(a, a_1; &\L, \Om)\bigr) 
	= \a^d \bigl(\GW_0(g(\re a); \L, \Om) +
		\GW_0(g(\re a_1); \L, \Om_1)
		\bigr) \\
&\ + \a^{d-1} \log\a \ \GW_1\bigl(
		\GD(g; \re a, \re a_1); \p\L, \p\Om \bigr) + o(\a^{d-1}\log\a),\ \a\to\infty.
	\end{align*}
The derivation of this 
fact from \eqref{jump:eq} repeats almost word for word the proof of formula \eqref{Widom1:eq}. 

Similarly to Theorem \ref{Szego3:thm}, Theorem \ref{Szego4:thm} 
is derived from formula \eqref{jump:eq} for polynomials $g_p(t) = t^p$, 
$p= 1, 2, \dots$:

\begin{thm}\label{jump_poly:thm}
Let the domains $\L, \Om$, and the symbols $a, a_1$ satisfy the 
conditions of Theorem \ref{Szego4:thm}. Then for any $p  = 1, 2, \dots$ we 
have 
\begin{align}\label{jump_poly:eq}
\lim_{\a\to\infty}\frac{1}{\a^{d-1}\log\a}
\tr\bigl[
&\ g_p\bigl(V_\a(a, a_1; \L, \Om)\bigr) 
- \chi_\L g_p\bigl( V_\a(a, a_1; \Rd, \Om)\bigr)\chi_\L
\bigr]\notag\\[0.2cm]
&\ \quad\quad  =  \GW_1\bigl(
\GD(g_p; a, a_1); \p\L, \p\Om \bigr).
\end{align}	
If $V_\a$ is replaced by  $H_\a$, then the same formula \eqref{jump_poly:eq} 
holds with $a, a_1$ replaced with $\re a, \re a_1$.
\end{thm}

The derivation of Theorem \ref{Szego4:thm} from Theorem \ref{jump_poly:thm} 
follows the plan of the proof of Theorem \ref{Szego3:thm}, and is omitted. 

As far as Theorem \ref{jump_poly:thm} itself is concerned, 
its proof essentially repeats that of Proposition \ref{Widom:prop} 
(given in \cite{Sob2}) with some 
obvious modifications. Below we provide only a sketch of this proof, 
leaving out some details that can be easily reconstructed.  

We'll need the following estimates.

%
%

\begin{prop}\label{cross1:prop}
Let the domain $\L$ be as in Theorem 
\ref{Szego4:thm}, and let 
the symbol $a\in\plainC\infty_0(\Rd\times\Rd)$  
be compactly supported in both variables. 
Then for any $q\in (0, 1]$, and all $\a\ge 1$ we have 
\begin{equation*}
\|\chi_\L \op_\a(a)(I-\chi_\L)\|_{\GS_q}^q \le C_q \a^{d-1},
\end{equation*}
and
\begin{equation*}
	\|P_{\Om, \a} \op_\a(a)(I-P_{\Om, \a})\|_{\GS_q}^q \le C_q \a^{d-1},
\end{equation*}
with a constant independent of $\a$.
\end{prop}

See \cite{Sob1}, Corollary 4.4.

\begin{proof}[Proof of Theorem \ref{jump_poly:thm} (sketch)] 
We give the proof only for the operator $V_\a$. The version 
of \eqref{jump_poly:eq} for the self-adjoint 
operator  $H_\a$ can be obtained following the elementary 
argument detailed in \cite{Sob}, p. 77. 
Furthermore, 
for simplicity we assume that $\L$ is a bounded domain, so that 	 
\eqref{jump_poly:eq} amounts to 
\begin{align}\label{main_inf_h:eq}
\lim_{\a\to\infty}\frac{1}{\a^{d-1}\log\a}
	 		\bigl[
	 		\tr \bigl(g_p(V_\a(a, a_1; \L, \Om))\bigr) 
	 		-  \a^d \bigl(
	 		\GW_0( &\ g_p(a_1); \L, \Om) + \GW_0(g_p(a_1); \L, \Om_1)\bigr)
	 		\bigr]\notag\\[0.2cm]
	 		=  &\ \GW_1(\GD( g_p; a, a_1); \p\L, \p\Om),
	 	\end{align} 
see the remark after Theorem \ref{Szego4:thm}. 
Since both $\L$ and $\Om$ are bounded, we may assume that the symbols $a, a_1$ 
are compactly supported in both variables. 
In what follows we use the following convention: for any 
two operators $A, B$ depending 
on the parameter $\a \ge 1$ 
we write $A\thicksim B$ if $\|A-B\|_{\GS_1}\le C\a^{d-1}$ with a constant $C$ 
independent of $\a$. 

By Proposition \ref{cross1:prop}, 
\begin{equation*}
P_{\Om, \a} a P_{\Om, \a} + P_{\Om_1, \a} a_1 P_{\Om_1, \a}
\thicksim 
a P_{\Om, \a} + a_1 P_{\Om_1, \a}
= 		b P_{\Om, \a} + a_1,\ b = a - a_1. 
\end{equation*}
Expanding $g_p\bigl(V(a, a_1; \L, \Om)\bigr)$ and repeatedly 
using Proposition 
\ref{cross1:prop} again, we obtain that 
		\begin{equation*}
			g_p\bigl(V(a, a_1; \L, \Om)\bigr)
			\thicksim  \sum_{l=0}^p {p \choose l} 
			\op_\a (a_1^{p-l}) g_l\bigl(
			T_\a(b; \L, \Om)
			\bigr).
		\end{equation*}
Traces of operators similar to the ones in the sum above 
have been studied in \cite{Sob2}. 
By \cite{Sob2} (see Lemma 3.3 and Theorem 4.1), 
	\begin{align}\label{nouveau:eq}
		\tr g_p\bigl(V(a, a_1; \L, \Om)\bigr)
		= &\ 
		\a^d \GW_0(a_1^p; \L, \Rd) + 
		\a^d\sum_{l=1}^p {p \choose l} \GW_0(a_1^{p-l}g_l(b); \L, \Om) 
		\notag \\[0.2cm]
		+ &\ \a^{d-1} \log \a \sum_{l=1}^p {p \choose l} 
		\GW_1\bigl(a_1^{p-l}\GA(g_l; b);\ \p\L, \p\Om\bigr)
		+ o(\a^{d-1}\log\a).
	\end{align}
	The second sum starts with $l=1$ since $\GA(g_0; b) = 0$. 
	By definition \eqref{w0:eq} 
	the first two terms on the right-hand side amount to 
	\begin{align*}
		\frac{\a^d}{(2\pi)^d}
		\int_\L\int_\Rd  &\ \bigl(b(\bx, \bxi) \chi_\Om(\bxi) + a_1(\bx, \bxi)
		\bigr)^p d\bxi d\bx \\[0.2cm]
		= &\ \a^d \bigl(\GW_0(g_p(a); \L, \Om) +
		\GW_0(g_p(a_1); \L, \Om_1)
		\bigr).
	\end{align*}
	To evaluate the second sum in \eqref{nouveau:eq} note that by definition \eqref{w1:eq},
	\begin{align*}
		\sum_{l=1}^p {p \choose l} z^{p-l}\GA(g_l; s)
		= &\ \frac{1}{(2\pi)^2}
		\int_0^1 \frac{g_p(z+ts) - g_p(z)- tg_p(z+s) + t g_p(z) }{t(1-t)} dt\\[0.2cm]
		= &\ \GD(g_p; s+z, z),
	\end{align*}  
	for any $z, s\in\mathbb C$. 
	Therefore the second sum on the right-hand side of \eqref{nouveau:eq} 
	equals 
	\begin{equation*}
		\a^{d-1} \log \a 
		\ \GW_1\bigl(
		\GD(g_p; a, a_1); \p\L, \p\Om\bigr).
	\end{equation*} 
	This completes the proof of \eqref{main_inf_h:eq}.  
	As explained earlier, this leads to 
	\eqref{jump_poly:eq}, as required. 
	\end{proof}

\bibliographystyle{amsplain}

\end{document}